# High-order adaptive multiresolution wavelet upwind schemes for hyperbolic conservation laws

Bing Yang, Jizeng Wang[*], Xiaojing Liu, Youhe Zhou

*Key Laboratory of Mechanics on Disaster and Environment in Western China (Lanzhou University), the Ministry of Education; College of Civil Engineering and Mechanics, Lanzhou University, Lanzhou, Gansu 730000, P.R. China*

**Abstract:** A system of high-order adaptive multiresolution wavelet collocation upwind schemes are developed for the solution of hyperbolic conservation laws. A couple of asymmetrical wavelet bases with interpolation property are built to realize the upwind property, and address the nonlinearity in the hyperbolic problems. An adaptive algorithm based on multiresolution analysis in wavelet theory is designed to capture moving shock waves and distinguish new localized steep regions. An integration average reconstruction method is proposed based on the Lebesgue differentiation theorem to suppress the Gibbs phenomenon. All these numerical techniques enable the wavelet collocation upwind scheme to provide a general framework for devising satisfactory adaptive wavelet upwind methods with high-order accuracy. Several benchmark tests for 1D hyperbolic problems are carried out to verify the accuracy and efficiency of the present wavelet schemes.

## 1  Introduction

Problems governed by hyperbolic conservation laws contain steep gradient regions or even discontinuities generally, such as shock waves and contact discontinuities, which is the core challenge in designing corresponding numerical schemes. In addition, the solutions of such problems often have substantially complex structures including smooth components such as vortices and acoustic waves. Traditional low-order numerical methods, such as the first-order Godunov scheme [1] and the Roe scheme [2], can capture discontinuities without spurious oscillations. However, they often smear the contact discontinuities excessively and carry comparatively large numerical dissipation in the smooth regions of the solutions. Therefore, a large number of nodes are required to recognize complex smooth structures accurately for long-time simulation, which is considerably inefficient.

High-order numerical schemes are more attractive in computational fluid dynamics because of their capability to achieve the desired resolution and higher accuracy with much fewer degrees of freedom [3]. Many types of high-order numerical methods for hyperbolic conservation laws have been investigated in the last decades due to their better potential in high efficiency. These schemes include finite difference method

[*] Corresponding Author. Email: jzwang@lzu.edu.cn

(FDM), finite volume method (FVM), finite element method (FEM), and meshless method. FDM is applied to problems on structured meshes over regular geometry because of simplicity and high efficiency. In unstructured meshes over complex geometry, FVM and FEM are more popular owing to their flexibility on arbitrary meshes.

Essentially non-oscillatory (ENO) schemes, first proposed in the frame of the FVM by Harten et al. [4] and developed in the FDM framework by Shu and Osher [5], fundamentally progress in obtaining high-order accuracy in smooth regions and sharp discontinuity transitions without spurious oscillations. The ENO schemes choose the smoothest stencil among several candidates to approximate the fluxes at the cell boundaries with high-order accuracy and to avoid spurious oscillations, leading to inefficiency and reduction of accuracy for certain functions. The weighted ENO (WENO) scheme in the FVM framework proposed by Liu et al. [6] is a reasonable remedy to overcome the shortcomings while maintaining the robustness and high-order accuracy of the ENO schemes. Jiang and Shu [7] extended the idea to the FDM and proposed the efficient implementation of WENO schemes, which greatly succeeded in high-order accuracy and high efficiency. Until now, the ENO and WENO schemes have been applied in many physical areas, such as turbulent flows [8], supersonic and hypersonic flows [9], astrophysical flows, atmospherical and climate sciences [10], and computational biology [11].

Another appealing type of the high-order methods is the compact scheme with localized data structures resulting in very high efficiency for parallel computing. Main examples of compact schemes are multimoment constrained finite volume (MCV) methods, spectral volume (SV) methods, and discontinuous Galerkin (DG) FEMs. The MCV schemes proposed by Ii and Xiao [12] define the point values within a single cell at equally spaced points and use a set of constraint conditions imposed on different kinds of moments to achieve the desired high-order accuracy. Wang et al. [13] devised the SV schemes by subdividing each basic volume into small control volumes and applying cell-averaged data from these control volumes to reconstruct a high-order approximation in the basic volume. The DG schemes in the FEM framework for time evolution problems proposed by Cockburn and Shu [14-16] convert the partial differential equations to a weak form by the conventional Galerkin method in one element and seek a discontinuous approximate solution in a space of polynomials. The DG methods borrow the Riemann solver of the classic FDM and FVM to obtain the numerical fluxes at the element boundaries [17]. The above compact schemes all apply the approximate Riemann solver to obtain the fluxes at the boundaries of cells or elements and can perform well on unstructured meshes over arbitrary geometries when dealing with hyperbolic conservation laws. However, the numerical solutions are still polluted by spurious oscillations near the discontinuities once the shocks appear inside the cells or elements. A large number of additional ingredients have been designed to remove these oscillations, for example, total variation bounded (TVB) limiters [12, 18] and WENO limiters [10, 19, 20].

The classical Galerkin FEM is ineffective when flows involve shocks or sharp layers [21]. Hughes et al. [22-24] developed streamline-upwind/Petrov Galerkin (SUPG) schemes with an extra discontinuity capturing term to advance the FEM progress in resolving hyperbolic problems. Unlike the DG schemes using the Riemann solver, the SUPG schemes achieve the upwind property based on modifying the weighted function with an additional term defined by the scalar product of velocity and the weighted function gradient [22]. The SUPG schemes stimulate the birth of other schemes in the FEM framework, such as Galerkin/least-square methods [25] and meshless local Petrov Galerkin method [26].

Multiscale smooth structures and discontinuities often appear simultaneously in the fluid field governed by hyperbolic conservation laws. The above presented high-order schemes are generally designed on uniform nodes or cells. Thus, a large number of nodes or cells must be set to detect and capture the different scales accurately, which is not cost effective. Adaptive algorithms with fine resolution near the localized steep regions and coarser resolution elsewhere provide a good ingredient to achieve high accuracy by reducing computational cost. Motivated by this goal, a large number of adaptive methods have been developed and applied to resolve the partial differential equations, for example, adaptive mesh refinement [27], FEMs [28, 29], and adaptive wavelet methods [30-35]. Posteriori or heuristic approaches as the criteria for controlling mesh refinement are generally used in the traditional adaptive methods. This approach would be quite costly for hyperbolic and parabolic problems because an iterative process is required to ensure the characteristic waves inside the refined area during one time step [36].

Multiresolution analysis in wavelet theory provides an extremely effective, common approach to capture the local characteristic of functions in time and frequency domains. Moreover, many wavelet-based numerical methods have been used in solving different types of PDEs [37, 38]. However, comparatively few works are dedicated for hyperbolic PDEs, which can be categorized into two types. The first one only embeds the wavelet multiresolution analysis into the traditional high-order schemes, for example, the DG schemes [39], the finite difference WENO [40], and the FDM with artificial viscosity [35]. The primary thought of the mixed methods is to detect the trouble nodes or cells near the localized steep structures and implement the adaptive or reconstruction process at the trouble ones based on wavelet coefficients. The remainder is called a pure wavelet numerical method that uses the wavelets as a set of bases to discretize the PDEs. Juan and Gray [41] first applied the classic wavelet Galerkin schemes in uniform node distribution to solve hyperbolic equations and found that spurious oscillations would spread further away from the shock, leading to numerical instability, which was also encountered in other classic Galerkin schemes. Recently, a dynamical wavelet Galerkin scheme was proposed by Pereira et al. [42] that overcomes the above drawback successfully due to the energy dissipation introduced by a non-smooth projection operator. Their work opens perspectives for studies of nonlinear hyperbolic conservation laws using adaptive Galerkin discretizations. In addition, Minbashian et al. [43] developed an adaptive wavelet SUPG method for hyperbolic conservation laws and performed a post processing using a denoising technique based on a minimization formulation to reduce Gibbs oscillations. Compared with the wavelet Galerkin methods, the collocation ones are more efficient. Adaptive wavelet collocation methods are applied to solve parabolic problems successfully based on symmetrical wavelet basis [44-46]. However, no attempts are carried out to use the adaptive wavelet collocation method for the hyperbolic problems. The basic theory of the hyperbolic PDEs show that a natural, fundamental thought is to design upwind schemes when solving the hyperbolic problems. The upwind schemes introduce an appropriate implicit numerical viscosity to ensure the numerical solution stably converges to the entropy solution [47]. Inspired by this point, a high-order adaptive multiresolution wavelet collocation upwind scheme for hyperbolic conservation laws that incorporates the merits of high-order schemes and adaptive algorithm is proposed in this paper.

The rest of the paper is organized as follows. Wavelet approximation theory, wavelet collocation upwind schemes for uniform and adaptive node distributions, and oscillation limiters are introduced in Section 2. Numerical experiments for the linear scalar equation, the nonlinear inviscid Burger's equation, and the Euler system in one dimension are carried out based on the proposed wavelet collocation upwind schemes in Section 3. The main conclusions are drawn in Section 4.

## 2 Fundamental theory and algorithms

One-dimensional scalar conservation laws are considered as a start point to construct wavelet upwind schemes:

$$u_t + f(u)_x = 0. \tag{2.1}$$

Before presenting the wavelet upwind schemes, wavelet approximation theory as preliminary knowledge is elaborated. Then, the basic ideas and formulations of the schemes are discussed in detail.

2.1 Wavelet approximation theory

2.1.1 Preliminaries

The Hilbert space $\mathbf{L}^2(\mathbf{\Omega})$ is the collection of square integrable functions that satisfy $\int_\Omega |f(x)|^2 \, \mathrm{d}x < \infty$, where $\mathbf{\Omega}$ is any open subset of $\mathbf{R}$ [48]. The space is naturally equipped with the inner product

$$(f, g) := \int_\Omega f(x) g(x) \, \mathrm{d}x < \infty, \quad \text{for all } f, g \in \mathbf{L}^2(\mathbf{\Omega}). \tag{2.2}$$

The associated $\mathbf{L}^2(\mathbf{\Omega})$ norm is defined by

$$\|f\|_{\mathbf{L}^2(\Omega)} = \left( \int_\Omega |f(x)|^2 \, \mathrm{d}x \right)^{1/2}, \quad \text{for all } f \in \mathbf{L}^2(\mathbf{\Omega}). \tag{2.3}$$

The $\mathbf{L}^\infty(\mathbf{\Omega})$ norm is also defined by

$$\|f\|_{\mathbf{L}^\infty(\Omega)} = \sup_{x \in \Omega} \{f(x)\}. \tag{2.4}$$

For any integer $m \geq 1$, the space of all functions that are $m$ times continuously differentiable over $\mathbf{\Omega}$ is denoted by $\mathbf{C}^m(\mathbf{\Omega})$.

To characterize the singular structures of functions, Lipschitz exponent is a good choice that provides point-wise or uniform regularity measurement over $\mathbf{\Omega}$. A function $f$ is Lipschitz $\alpha$ at the point $x$ if it satisfies

$$|f(x) - f(y)| \leq \mathrm{K} |x - y|^\alpha \quad \text{for } x \in \mathbf{\Omega}, \text{ any } y \in (x - \delta, x + \delta), \tag{2.5}$$

where K is a constant that is dependent on the point $x$, and $\delta > 0$ in a small magnitude [49]. The present paper concentrates on the wavelet transform defined on $\mathbf{R}$ and functions in $\mathbf{L}^2(\mathbf{\Omega}) \cap \mathbf{C}^m(\mathbf{\Omega})$.

Our attention now turns to the primary theory of wavelets beginning with the definition of biorthogonal multiresolution analysis. A biorthogonal multiresolution analysis consists of two dual sequences of closed subspaces of $\mathbf{L}^2$ denoted by $V_J$ and $\tilde{V}_J$, which satisfy the relations

$$\begin{aligned} \cdots V_{-1} \subset V_0 \subset V_1 \subset \cdots \subset V_J \subset V_{J+1} \subset \cdots, \\ \cdots \tilde{V}_{-1} \subset \tilde{V}_0 \subset \tilde{V}_1 \subset \cdots \subset \tilde{V}_J \subset \tilde{V}_{J+1} \subset \cdots, \end{aligned} \tag{2.6}$$

$$\overline{\bigcup_{J \in \mathbf{Z}} V_J} = \overline{\bigcup_{J \in \mathbf{Z}} \tilde{V}_J} = \mathbf{L}^2, \tag{2.7}$$

$$\overline{\bigcap_{J \in \mathbf{Z}} V_J} = \overline{\bigcap_{J \in \mathbf{Z}} \tilde{V}_J} = \{0\}, \tag{2.8}$$

$$f \in V_J \Leftrightarrow f\left(2^{-J} \cdot\right) \in V_0,$$
$$f \in \tilde{V}_J \Leftrightarrow f\left(2^{-J} \cdot\right) \in \tilde{V}_0,$$
(2.9)

$$f \in V_0 \Leftrightarrow f(\cdot - k) \in V_0,$$
$$f \in \tilde{V}_0 \Leftrightarrow f(\cdot - k) \in \tilde{V}_0,$$
(2.10)

where *J* is the resolution level. Generally, two different auxiliary functions, $\varphi$ and $\tilde{\varphi}$, known as the scaling function and its corresponding dual scaling function, respectively, satisfy

$$\{\varphi(\cdot - k)\}_k \text{ is a Riesz basis of } V_0,$$
$$\{\tilde{\varphi}(\cdot - k)\}_k \text{ is a Riesz basis of } \tilde{V}_0,$$
$$(\varphi(\cdot - l), \tilde{\varphi}(\cdot - k)) = \delta_{l,k}.$$
(2.11)

The filter coefficients of the scaling functions can be determined based on refinement relation

$$\varphi(x) = \sum_k h_k \varphi(2x - k),$$
$$\tilde{\varphi}(x) = \sum_k \tilde{h}_k \tilde{\varphi}(2x - k).$$
(2.12)

The scaling functions at level *J* can be defined by

$$\varphi_{J,k}(x) = 2^{J/2} \varphi\left(2^J x - k\right),$$
$$\tilde{\varphi}_{J,k}(x) = 2^{J/2} \tilde{\varphi}\left(2^J x - k\right).$$
(2.13)

Following the above scaling functions, the wavelets are defined as

$$\psi(x) = \sum_k g_k \varphi(2x - k),$$
$$\tilde{\psi}(x) = \sum_k \tilde{g}_k \tilde{\varphi}(2x - k),$$
(2.14)

where

$$g_k = (-1)^k \tilde{h}_{1-k},$$
$$\tilde{g}_k = (-1)^k h_{1-k}.$$
(2.15)

Similarly, the wavelets at level *J* are given by

$$\psi_{J,k}(x) = 2^{J/2} \psi\left(2^J x - k\right),$$
$$\tilde{\psi}_{J,k}(x) = 2^{J/2} \tilde{\psi}\left(2^J x - k\right).$$
(2.16)

The spaces spanned by $\psi_{J,k}(x)$ and $\tilde{\psi}_{J,k}(x)$ are denoted by $W_J$ and $\tilde{W}_J$, respectively. The wavelet spaces represent the difference between two successive approximations, $V_{J+1}$ and $V_J$, that means

$$V_{J+1} = V_J \oplus W_J,$$
$$\tilde{V}_{J+1} = \tilde{V}_J \oplus \tilde{W}_J. \tag{2.17}$$

All functions in $\mathbf{L}^2(\Omega)$ can be approximated arbitrarily close by the scaling functions because $V_J$ and $\tilde{V}_J$ are nested and dense in $\mathbf{L}^2$. Then, a projection operator $P_J : \mathbf{L}^2 \to V_J$ is found so that for every $f$ in $\mathbf{L}^2(\Omega)$ $\|P_J f - f\|_{\mathbf{L}^2}$ tends to 0 as $J$ approaches infinity. Next, the projection is defined in the biorthogonal setting as

$$P_J f = \sum_k (f, \tilde{\varphi}_{J,k}) \varphi_{J,k}. \tag{2.18}$$

Given $V_J = V_{j_0} \oplus \bigoplus_{j=j_0}^{J-1} W_j$, $P_J f$ can also be written in a multiresolution decomposition form as

$$P_J f = \sum_k c_{j_0,k} \varphi_{j_0,k} + \sum_{j=j_0}^{J-1} \sum_k d_{j,k} \psi_{j,k}, \tag{2.19}$$

where $c_{j_0,k} = (f, \tilde{\varphi}_{j_0,k})$, $d_{j,k} = (f, \tilde{\psi}_{j,k})$.

2.1.2 Wavelet construction

So far, the multiresolution analysis and function expansion in wavelet spaces have been presented. Next, how to construct a desirable wavelet used in a particular application is introduced. This work is limited to the construction of interpolation wavelets that belong to the family of second-generation wavelets. For more general wavelets, the detailed building methods can be found in the references [49-52].

Owing to our interest in numerical calculation and analysis, the scaling functions of "good" wavelets should have the following properties:

(1) Compact support: $\varphi(x)$ is only nonzero in a finite interval [$N_L$, $N_R$], which corresponds to the locality of the subdivision scheme.

(2) Interpolation: $\varphi(x)$ is interpolation in the sense that $\varphi(k) = \delta_{0,k}$, where $\delta_{i,j}$ is the Kronecker delta function.

(3) Polynomial reproduction: The scaling functions can reproduce polynomials up to the degree $N - 1$.

(4) Smoothness: $\varphi(x) \in \mathbf{C}^\alpha$, where $\alpha = \alpha(N)$. The smoothness increases linearly with $N$ [53].

(5) Refinability: The refinement relation is referred to as Equation (2.12).

For general applications, an additional property is symmetry, especially in signal processing and data compression. Several techniques are used to devise scaling functions satisfying the above conditions, for example, the auto correlation of Daubechies scaling functions [54, 55] and interpolation methods [56]. The auto correlation method can only be applied to build symmetric scaling functions. The latter provides more freedom to construct symmetric and asymmetric scaling functions. Then, a step-by-step procedure of building wavelets based on the interpolation methods is discussed:

(1) Choose a type of wavelet and determine the smoothness parameter $N$: (a) symmetrical wavelet, $N \in$ *Even*

is the only choice; (b) asymmetrical wavelet, $N \in odd$ or $N \in Even$ is optional.

(2) Select a node stencil. Specify $N$ nodes uniformly distributed in $V_J$ as the base and choose one node in $W_J$. The stencil for the symmetrical wavelet is unique, and several candidates are allowed for the asymmetrical one. To show stencils intuitively, $N = 6$ and $N = 5$ are chosen for examples:

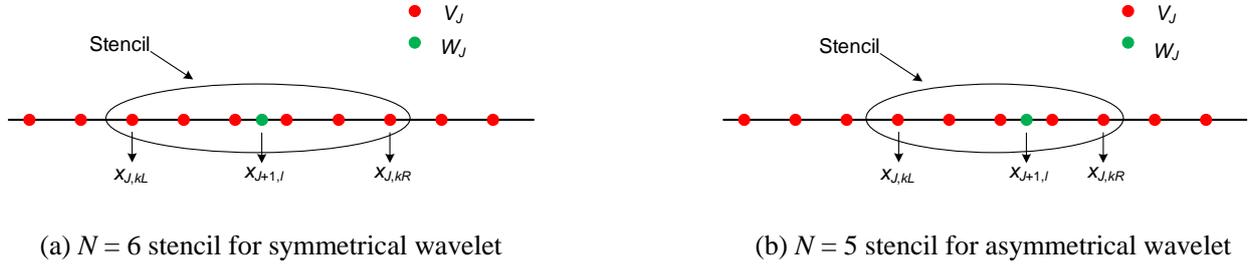

(a) $N = 6$ stencil for symmetrical wavelet  (b) $N = 5$ stencil for asymmetrical wavelet

Fig.1. Example for stencils

(3) Approximate $\varphi_{J,k}(x_{J+1,l})$ by $(N-1)$th order Lagrange interpolation polynomial and calculate the filter coefficients. The following relation is derived:

$$\varphi_{J,k}(x_{J,k_1}) = \delta_{k,k_1},$$
$$\varphi_{J,k}(x_{J+1,l}) = \sum_{k_1=kL}^{kR} \varphi_{J,k}(x_{J,k_1}) L_{J,k_1}(x_{J+1,l}), \quad (2.20)$$

where $x_{J,k_1} = k_1/2^J$ and $x_{J+1,l} = l/2^{J+1}$. Based on the refinement relation,

$$\varphi_{J,k}(x_{J+1,l}) = \sum_{l_1} h_{l_1} \varphi_{J+1,l_1}(x_{J+1,l}). \quad (2.21)$$

Substituting (2.21) into (2.20), $h_l = \delta_{0,l}$ can be easily calculated for $l \in Even$, for $l \in Odd$ by

$$h_l = L_{J,k}(x_{J+1,l}) = \prod_{\substack{i=kL \\ i \ne k}}^{kR} \frac{x_{J+1,l} - x_{J,i}}{x_{J,k} - x_{J,i}}. \quad (2.22)$$

For a specified parameter $k$, parameters $kL$ and $kR$ vary with $l$ translating. To compute the coefficients more efficiently, supposing that $J = 0$ and $k = 0$ can obtain

$$h_l = \prod_{\substack{i=kL \\ i \ne 0}}^{kR} \frac{l/2 - i}{-i}. \quad (2.23)$$

(4) Calculate the derivatives and integrals of the scaling function by algorithms proposed by Wang [57] and Chen et al. [58]. Returning to the refinement relation and applying the cascade algorithm, the values of the scaling functions and its derivatives and integrals at dyadic points at arbitrary refined resolution levels can be obtained.

Following the above steps, a desirable scaling function can be successfully built. Then, the dual scaling function needs to be determined. For the present interpolation method, the Kronecker delta function is

naturally chosen as the dual scaling function [56]:

$$\tilde{\varphi}_{J,k} = \delta(x - x_{J,k}). \tag{2.24}$$

Next, the thought proposed by Donoho [54], which builds a wavelet function from the scaling function, is followed:

$$\psi(x) = \varphi(2x - 1),$$
$$\psi_{J,k}(x) = \varphi\left(2^{J+1}x - (2k+1)\right). \tag{2.25}$$

Wavelet coefficients should be calculated by $(f, \tilde{\psi}_{J,k})$, but here for simplicity, the wavelet coefficients are calculated by using multiresolution analysis:

$$d_{J,l} = c_{J+1,2l+1} - \sum_k h_{(2l+1)-2k} c_{J,k}. \tag{2.26}$$

2.1.3 Approximation of functions on a finite domain

The decomposition (2.19) of the function is defined on the whole real line **R**. For general practical cases, $\Omega$ is a finite domain. Then, a new trouble occurs near the boundaries [54, 55, 57]. A boundary extension technique was developed in our previous study based on Lagrange interpolation to remove the local errors induced by a loss of information outside the domain [55]. The modified approximation in level $J$ can be written as

$$P_J f(x) = \sum_{k \in \mathbb{R}_J} f(x_k) \Phi_{J,k}(x), \tag{2.27}$$

where the modified wavelet basis is given by

$$\Phi_{J,k}(x) = \varphi_{J,k}(x) + \sum_{n \in \bar{\Xi}_k} \Pi_k^n(x_n) \varphi_{J,n}(x)$$
$$= \sum_{n \in \Xi_k} \Pi_k^n(x_n) \varphi_{J,n}(x), \tag{2.28}$$

$$\Pi_k^n(x_n) = \prod_{\substack{i=1 \\ i \neq n}}^{\eta} \frac{x_n - \tilde{x}_i}{\tilde{x}_k - \tilde{x}_i}, \tag{2.29}$$

where $\bar{\Xi}_k$ is the set of serial numbers $n$ denoting the external point $x_n$ such that $f(x_n)$ is exactly all the Lagrange interpolations. In Equation (2.29) the set { $\tilde{x}_i$, $i=1, 2,…,\eta$} is the collection of selected Lagrange interpolation nodes inside the $\Omega$ for the external point $x_n$ and $\Xi_k = \bar{\Xi}_k \cup \{k\}$ with the coefficient $\Pi_k^k(x_k) \equiv 1$. The modified wavelet basis $\Phi_{J,k}(x)$ keeps the interpolation property because the original scaling function and the Lagrange interpolation polynomial have such characteristics. Next, a series of wavelet approximation theorems without proof, which can be proven by using the approach proposed in the references, are presented [55, 59].

following relations:

$$\left\| P_{J_0}^{J_{\max}} f - f \right\|_{\mathbf{L}^2(\Omega)} \leq C_{4,L} 2^{-J_0 \lambda}, \quad (2.37)$$

$$\left\| P_{J_0}^{J_{\max}} f - f \right\|_{\mathbf{L}^\infty(\Omega)} \leq C_{4,\infty} 2^{-J_0 \lambda}. \quad (2.38)$$

Remarks: (1) All the function values $f(x_{J,k})$ used in the multiresolution approximation (2.36) satisfy the identical relation $P_{J_0}^{J_{\max}} f(x_{J,k}) \equiv f(x_{J,k})$. (2) The boundary extension is conducted at the basic level $J_0$, and $\eta = N$ is recommended as an appropriate choice. (3) Additional nodes at high resolution levels can be completely random. (4) **Theorem 2.4** shows that the additional nodes would never reduce the fundamental accuracy approximated based on the uniform distribution of nodes.

2.2 Wavelet collocation upwind schemes for uniform node distribution

This part shows how to establish a wavelet collocation upwind scheme with a uniform node distribution denoted by the WCU. The concept of upwind scheme in the FDM, FVM, and DG schemes is essential when devising new numerical schemes in solving hyperbolic conservation systems [7, 12-14, 60-62]. Here the Galerkin method is dismissed, and a first attempt is made to construct a pure wavelet collocation upwind scheme to resolve the hyperbolic conservation laws.

An asymmetrical wavelet may have upwind property. When calculating the derivative of the function at a discrete point based on differential operator of the wavelet approximation, the derivatives are obtained by a linear combination of several function values $u_{J,k}$:

$$u'(x_{J,l}) = \sum_{k \in \mathbb{R}_J} u_{J,k} \Phi'_{J,k}(x_{J,l}), \quad (2.39)$$

where $\Phi'_{J,k}(x_{J,l})$ is the derivative of the scaling function $\Phi_{J,k}(x)$ at $x_{J,l}$. The summation in Eq. (2.39) is composed of finite terms due to the compact support of the scaling functions. Inspired by the upwind scheme in classic numerical methods, a similar upwind scheme may be established as if choosing an appropriate asymmetrical wavelet. Asymmetry of wavelet means that the derivative of the scaling function is non-identical on two sides of the zero point and shows biased characteristics. To construct a high-order, stable scheme, the key point is to build a couple of "good" asymmetry wavelets.

The upwind property is determined by the location of the point chosen in $W_J$. If the node number of the stencil on the left side of $x_{J+1,l}$ is more than that on the right side, the wavelet is upwind in the positive direction. Otherwise, it is upwind in the negative direction. Two identically mirror symmetrical wavelets can be built to construct the upwind scheme in the positive and negative directions, which is similar to the upwind scheme in the FDM. Next, how to build a couple of wavelets for a specified parameter $N$ is shown.

Fig. 2 shows a stencil in $V_0$. The nodes in the stencil are symmetrical with zero point. When $x_{1,1}$ is placed between 0 and 1, a slightly asymmetrical wavelet with positive upwind property is obtained. The corresponding negative upwind wavelet can be easily constructed by mirroring the point in $W_0$ with zero point. By this strategy, the filter coefficients are also of mirror symmetry. This approach does improve our efficiency because one group of coefficients needs to be calculated. Here nonzero filter coefficients are shown in Table 1, and the scaling functions of two couples of wavelets are given as examples, as shown in

Fig. 3.

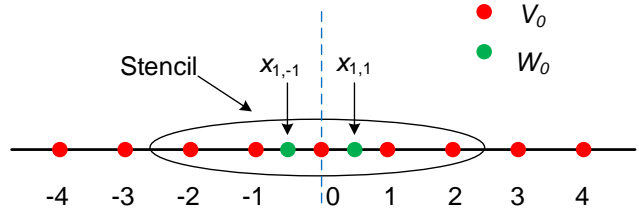

Fig.2. Stencils for a couple of asymmetrical wavelets

Table 1. Nonzero filter coefficients of some wavelets

| $N = 5$ | | $N = 7$ | |
|---|---|---|---|
| Positive upwind | Negative upwind | Positive upwind | Negative upwind |
| $h_{-3} = -0.0390625000$ | $h_{-5} = 0.023437500$ | $h_{-5} = 0.0068359375$ | $h_{-7} = -0.0048828125$ |
| $h_{-1} = 0.4687500000$ | $h_{-3} = -0.1562500000$ | $h_{-3} = -0.0683593750$ | $h_{-5} = 0.0410156250$ |
| $h_0 = 1.0000000000$ | $h_{-1} = 0.7031250000$ | $h_{-1} = 0.5126953125$ | $h_{-3} = -0.1708984375$ |
| $h_1 = 0.7031250000$ | $h_0 = 1.0000000000$ | $h_0 = 1.0000000000$ | $h_{-1} = 0.6835937500$ |
| $h_3 = -0.1562500000$ | $h_1 = 0.4687500000$ | $h_1 = 0.6835937500$ | $h_0 = 1.0000000000$ |
| $h_5 = 0.02343750000$ | $h_3 = -0.0390625000$ | $h_3 = -0.1708984375$ | $h_1 = 0.5126953125$ |
| | | $h_5 = 0.0410156250$ | $h_3 = -0.0683593750$ |
| | | $h_7 = -0.0048828125$ | $h_5 = 0.0068359375$ |

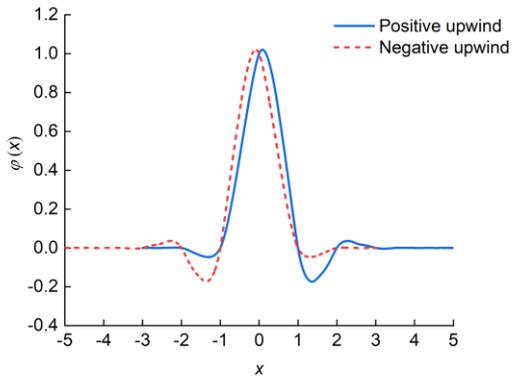

(a) $N=5$ Scaling function

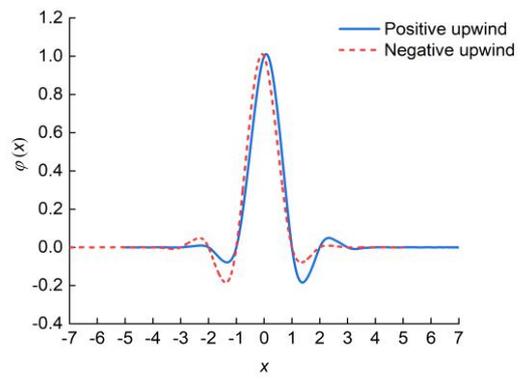

(b) $N=7$ Scaling function

Fig.3. Scaling functions of the two couples of wavelets

For problems governed by the hyperbolic conservation laws, the FDM suggests to solve these problems in the conservative form as shown in Equation (2.1). Fortunately, the nonlinear terms in the conservative equations can easily be discretized because of the exact interpolation property of the scaling function, which enables performing the following operator [63]:

$$f(u(x)) = \sum_{k \in \mathbb{R}_J} f(u_{J,k}) \Phi_{J,k}(x), \tag{2.40}$$

where $f$ satisfies a uniform Lipschitz condition of order $\alpha$ with respect to $u$, $\alpha \geq 1$. The decomposition coefficients of $f(u)$ can be evaluated by computing value of $f(u_{J,k})$ in a quite high efficiency.

Let us return to Equation (2.1). Considering a periodic boundary condition and $f'(u) \geq 0$, space discretization is carried out at the basic resolution level $J$ by the wavelet collocation upwind scheme, and the following semi-discretized system is obtained:

$$\frac{du_J(x_l, t)}{dt} + \sum_{k \in \mathbb{R}_J} f(u_J(x_k, t)) \Phi'_{+,J,k}(x_l) = 0, \quad k, l \in \Omega, \tag{2.41}$$

where $\Phi_{+,J,k}(x)$ denotes the scaling function of the positive upwind wavelet. Then, the above ordinary partial equations can be solved by a general method. In this paper, the classic explicit fourth-order Runge–Kutta method conducting time integration is applied. In the frame of the uniform nodes, the derivative $\Phi'_{J,k}(x_l)$ of the scaling function is a series of unchanged coefficients for the specified level $J$. This result indicates that our schemes are conservative.

Now our schemes are applied to a general flux $f(u)$, that is, $f'(u) \ngeq 0$ that can be split into two parts either globally or locally:

$$f(u) = f^+(u) + f^-(u), \tag{2.42}$$

where $df^+(u)/du \geq 0$ and $df^-(u)/du \leq 0$. The flux splitting techniques developed in FDM to realize the space discretization for the Euler system is further borrowed because our schemes follow the upwind concept. Several splitting methods are successful, such as Steger warming flux splitting [64], Lax–Friedrichs flux splitting [7], and Roe flux splitting [65]. Here Global Lax–Friedrichs flux splitting is chosen and can be described as follows:

$$f^{\pm}(u) = \frac{1}{2}(f(u) \pm \alpha u), \tag{2.43}$$

where $\alpha = \max |f'(u)|$, and the maximum is taken over the whole relevant range of $u$. Then, the general semi-discretized system can be obtained:

$$\frac{du_J(x_l, t)}{dt} + \sum_{k \in \mathbb{R}_J} f^+(u_J(x_k, t)) \Phi'_{+,J,k}(x_l) + \sum_{k \in \mathbb{R}_J} f^-(u_J(x_k, t)) \Phi'_{-,J,k}(x_l) = 0, \; k, l \in \Omega, \tag{2.44}$$

where $\Phi_{-,J,k}(x)$ is the scaling function of the negative upwind wavelet. To extend the scheme to 1D Euler system, a vector consisting of three conservative variables only needs to be substituted into Equation (2.44).

Based on the above analysis, our wavelet collocation upwind schemes solve the conservative form equations in a conservative way. Therefore, our schemes can work well in solving the hyperbolic conservative laws.

## 2.3 Adaptive multiresolution wavelet collocation upwind schemes

In this section, the detailed instruction on implementation of the proposed wavelet collocation upwind scheme is described. For simplicity, 1D scalar conservation laws are used to clarify the adaptive multiresolution wavelet collocation upwind scheme abbreviated as AMWCU. According to the general semi-discretized system in uniform nodes, the general semi-discretized system for the adaptive algorithm can be obtained:

$$\frac{du_J(x_l,t)}{dt} + \left( \sum_{k \in \mathbb{R}_{J_0}} f^+(u_{J_0}(x_k,t)) \Phi'_{+,J,k}(x_l) + \sum_{J=J_0+1}^{J_{max}} \sum_{k \in \mathbb{R}_J} d^+_{J,k} \Phi'_{+,J,k}(x_l) \right) \\ + \left( \sum_{k \in \mathbb{R}_{J_0}} f^-(u_{J_0}(x_k,t)) \Phi'_{-,J,k}(x_l) + \sum_{J=J_0+1}^{J_{max}} \sum_{k \in \mathbb{R}_J} d^-_{J,k} \Phi'_{-,J,k}(x_l) \right) = 0, \ k,l \in \Omega \quad (2.45)$$

where $d^+_{J,k}$ and $d^-_{J,k}$ are the wavelet coefficients associated with high resolution level of $f^+(u)$ and $f^-(u)$ regarding the positive and negative upwind wavelets, respectively.

The most appealing feature of wavelet analysis for solving PDEs is the ability to evaluate the local regularity of the solution based on the wavelet coefficients. The core thought of the adaptive wavelet algorithm is to conduct a self-adaptive discretization with automatic local node refinement and obtain a sparse representation of the unsolved variables. Returning to the multiresolution decomposition of the function (2.36), only the wavelet coefficients $|d_{J,k}| > \varepsilon$ are retained:

$$P_{J_0}^{J_{max}} u(x) = \sum_{k \in \mathbb{R}_{J_0}} u(x_k) \Phi_{J,k}(x) + \sum_{J=J_0+1}^{J_{max}} \sum_{\substack{k \in \mathbb{R}_J \\ |d_{J,k}| > \varepsilon}} d_{J,k} \Phi_{J,k}(x), \quad (2.46)$$

where $\varepsilon$ is the thresholding parameter with small magnitude, for example, $10^{-5}$. Here a node with wavelet coefficient $|d_{J,k}| > \varepsilon$ is named as a trouble node. Accurately approximating a function with local singularity by (2.46) is sufficient. However, solving hyperbolic evolution problems with shock discontinuities might be problematic. To ensure the adequate approximation of the solution during time evolution, Liandrat and Tchamitchian [66] proposed the notion of an adjacent zone including wavelets whose coefficients are or might become substantial during a step of time integration, when the node distribution remains unchanged in the step. Our node adaptive strategy is also based on this idea.

For general hyperbolic conservative laws, the domain $\Omega$ defining the problems can be divided into two parts. The first part has sufficiently smooth solution whose nodes are all at resolution level $J_0$ and another part with the steep gradient solution has nodes at high resolution levels. Shock might occur in such a smooth region, so methods to recognize the possible trouble nodes based on the information provided by $J_0$ resolution level must be devised. Therefore, two different adaptive strategies need to be designed, namely, adaptive node generation for $J_0$ (ANG-$J_0$) and adaptive node generation for $J$ (ANG-$J$).

First, the details of ANG-$J_0$ are shown. To measure the regularity of the solution in a smooth region, the smoothness measuring function proposed by Jiang and Shu is applied [7]:

$$IS_{J_0,l} = \frac{13}{12}\left(f_{J_0,l-1} - 2f_{J_0,l} + f_{J_0,l+1}\right)^2 + \frac{1}{4}\left(f_{J_0,l-1} - f_{J_0,l+1}\right)^2$$
$$= (f'\Delta x)^2\left(1 + O(\Delta x^2)\right) \quad \text{for } \Delta x = 2^{-J_0}, f' \neq 0. \tag{2.47}$$

One can see that $IS_{J_0,l}$ is affected by the space step $\Delta x$ and the first-order derivative of the function. Therefore, a suitable parameter $M_0$ should be determined to detect the location with low regularity, especially discontinuities. If $IS_{J_0,l} > M_0 \Delta x^2$, the node $x_{J_0,l}$ is marked as the trouble one. Once all trouble nodes are recognized, a set of nodes are inserted at the high resolution level based on the algorithm proposed in our previous research [67] to improve the local approximation accuracy. The strategy with suitable parameters adds enough nodes in the adjacent zone of the trouble nodes to guarantee the accuracy for time evolution problems.

Then, the ingredients for ANG-$J$ are presented. Based on the approach proposed by Vasilyev and Paolucci [45], a rule is developed to determine if neighboring nodes are in the adjacent zone of a trouble node. Supposing $x_{J,l}$ to be a trouble node, a node $x_{J_1,k}$ belonging to the adjacent zone of the node $x_{J,l}$ can be defined if the following relations are satisfied:

$$|J_1 - J| \leq L, |x_{J_1,k} - x_{J,l}| \leq K_w \cdot 2^{-J}, \tag{2.48}$$

where $L$ is the difference between the maximum resolution level of the added nodes and the level $J$, and $K_w$ is the width of the adjacent zone. The parameters $L$ and $K_w$ specify the total number and positions of the nodes that are inserted in the adjacent zone of the node $x_{J,l}$ at $t_i$, respectively. Here, $L = 1$ and $K_w = 2$ is recommended for most problems based on numerical tests. In addition, calculations are performed without modifying the node distribution when integrating the semi-discretization systems from $t_i$ to $t_{i+1}$.

After presenting the above fundamentals, the detailed procedures of implementing the dynamical process are described. Compared with the scheme of uniform node distribution, additional tasks in such an adaptive process are requirements of the dynamically adaptive node generation and calculation of wavelet coefficients of the fluxes. The process of dynamically adaptive node generation includes node initialization and dynamical adaption. The former has the following steps:

(1) Determine a characteristic variable that can depict all the singularity positions of the unsolved variables and calculate the wavelet coefficients. For 1D Euler system, the conservative variable $\rho$ is recommended.

(2) Choose an appropriate basic resolution level $J_0$ and maximum resolution level $J_{max}$, create uniform nodes for all $x_{J_0,k} \in \Omega$, and add two nodes at $x = a, b$, which are the endpoints of $\Omega$ to impose the boundary conditions.

(3) Detect the localized structure of the initial condition and insert nodes associated with high resolution level $J$ based on ANG-$J_0$. Then, delete the repeated nodes and external nodes, and obtain a new node distribution $\mathcal{N}^1$.

(4) Reconstruct the values of the conservative variables and carry out time integration to obtain $u(x_{J,l}, t_0 + \Delta t)$.

Then, adaptive node generation is performed at every time step:

(1) Obtain the conservative variables of the solution at $t_i$ ($i > 0$) with computational node distribution $\mathcal{N}^i$, choose a characteristic variable, and compute the wavelet coefficients.

(2) Evaluate $IS_{J_0,l}$ at the basic resolution level at $t_i$, and determine the total node number $N_{J_0}$ added in the new steep gradient area with ANG-$J_0$.

(3) Compare the wavelet coefficients at the high resolution level $J$ ($J_0 < J < J_{\max}$) with the specified threshold parameter $\varepsilon$, and gain the total node number $N_J$ in the adjacent zones of the trouble nodes applying ANG-$J$.

(4) Insert the added nodes based on ANG-$J_0$ and ANG-$J$ in $\mathcal{N}^i$, delete the repeated nodes and external nodes, and obtain an updated node distribution $\mathcal{N}^{i+1}$.

(5) Reconstruct the values of the conservative variables in $\mathcal{N}^{i+1}$ at $t_i$, and integrate the resulting semi-discretized system to the next step $t_{i+1}$. If $t_{i+1} = t_{\text{end}}$, then the solution is complete. Otherwise, repeat (1)–(5).

We now present an approach to calculate the wavelet coefficients in a fast, efficient way. Recalling the conventional fast wavelet transform algorithm [56], the forward transform process works level wise and on each level splits coefficients from the maximum level $J_{\max}$ to the basic level $J_0$. This algorithm does perform well when all information is given for the function. However, in numerical calculations, some unknown values of the functions need to be obtained by wavelet approximation when computing the wavelet coefficients based on multiresolution analysis. Our previous study [55] presented that the standard interpolation wavelet basis function demonstrates the property $\Phi_{J,k}(x_{J',l}) = 0$ for $J > J'$. This relation indicates that the wavelet coefficients at higher resolution level have no effect on the function values of the dyadic points at the lower resolution level. Based on this conclusion, only the decomposition direction of the forward transform is reversed, which means that the fast wavelet transform is conducted from the basic level $J_0$ to the maximum level $J_{\max}$. Then, the following fast forward transform is provided:

(1) Set the initial value: $u_{J_0,l} = u_{J_0,l}, d_{J,l} = u_{J,l}$ for $J > J_0$, and let $J = J_0 + 1$.

(2) Select the nodes $x_{J,l}$ in $J$, and compute $d_{J,l}$ based on the formula:

$$d_{J,l} = d_{J,l} - \sum_k h_{(2l+1)-2k} u_{J-1,k}. \tag{2.49}$$

If the node set excludes $x_{J-1,k}$ from the $\mathcal{N}^{i+1}$, calculate $u_{J-1,k}$ by the wavelet approximation.

(3) Compare $J$ with $J_{\max}$. If $J = J_{\max}$, then the process ends; otherwise, let $J = J + 1$ and repeat (2)–(3).

2.4 Gibbs phenomenon suppression

When approximating a function with a jump discontinuity using a set of smooth basis, numerical oscillations always exist near the discontinuity, which is similar to the Gibbs phenomenon in the Fourier series [68]. Therefore, the Gibbs phenomenon is inevitable when approximating functions with the jump discontinuity based on the wavelets. Thus, a remedy to remove the non-physical spurious oscillations when capturing shock waves in solving hyperbolic conservative laws must be proposed.

**Theorem 2.5 (Lebesgue differentiation theorem)** [69] If $f$ is Lebesgue integrable on $\mathbf{R}^d$, then

$$\lim_{\substack{m(B)\to 0 \\ x\in B}} \frac{1}{m(B)} \int_B f(y)dy = f(x) \text{ for a.e. } x. \tag{2.50}$$

The Lebesgue differentiation theorem indicates that the integral average value of the function over balls $B$ containing $x$ becomes $f(x)$ as $m(B)$ approaches 0. Enlightened by this theorem and TVB limiters in classical methods [18, 70], an integral reconstruction method is proposed to eliminate the spurious oscillations.

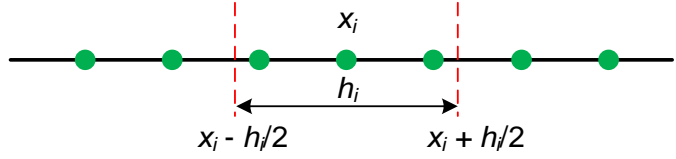

Fig.4. Sketch of the integral interval

Choosing an arbitrary node $x_i \in \Omega$ as shown in Fig. 4, the integral average value of $u$ over the interval $[x_i - h_i/2, x_i + h_i/2]$ is defined by

$$\bar{u}_i = \frac{1}{h_i} \int_{x_i-h_i/2}^{x_i+h_i/2} u(x)\, dx, \tag{2.51}$$

where $h_i$ is the length of the integral interval. By conducting Taylor expansion, the following can be obtained:

$$\bar{u}_i = u_i + O(h_i^2). \tag{2.52}$$

A switch function is defined as

$$\tilde{m}(a_1) = \begin{cases} -a_1 & \text{if } |a_1| \leq Mh_i^2 \\ 0 & \text{otherwise} \end{cases}, \tag{2.53}$$

where M is a switching parameter. Then, $u_i$ is reconstructed by the relation

$$u_i = \bar{u}_i + \tilde{m}_1(\bar{u}_i - u_i). \tag{2.54}$$

We can see that $u_i$ is only modified by the integral average value when $|\bar{u}_i - u_i| > Mh_i^2$. For smooth regions, the quantity of $|\bar{u}_i - u_i|$ is quite small, and $u_i$ maintains its high accuracy. The accuracy of the integral remarkably affects the reconstruction process. However, the wavelet approximation of integrals has high-order accuracy, as described in **Theorem 2.3,** which is a good choice for computing the integral.

We remark that $h_i$ is a crucial parameter and problem dependent. If $h_i$ is so large, the accuracy may be reduced in a smooth region. However oscillations are eliminated unsuccessfully when $h_i$ is too small. Moreover, $h_i$ should be determined by $K_1 \cdot 2^{-J}$ to ensure integral accuracy. The unsolved problems can be classified into two categories. First, the initial condition contains discontinuities, which permits refining nodes beforehand. Second, the initial distribution is sufficiently smooth, and the singularity of the solution emerges as time varies. As to the former case, $h_i = 2^{-J_{\max}}$ is recommended for scalar conservation laws and Euler systems ($J_{\max} = J_0$ in the uniform node algorithm). For the latter, determining an appropriate $h_i$ is slightly difficult. Wavelet coefficient is a natural discontinuity indicator, so a rule can be designed by comparing the value with predefined parameters to control the quantity of $h_i$. A reasonable option to compute the length of the integral interval is provided as follows:

$$h_i = \begin{cases} \dfrac{1}{2^{\lfloor 0.5(J_{\max}+J_0) \rfloor}} & J_i > J_0, d_{J_i,i} \geq \varepsilon_0, \\ 0.5\left(\dfrac{1}{2^{\lfloor 0.5(J_{\max}+J_0) \rfloor}} + \dfrac{1}{2^{J_{\max}-1}}\right) & J_i > J_0, \varepsilon_1 \leq d_{J_i,i} < \varepsilon_0, \\ \dfrac{1}{2^{J_{\max}-1}} & J_i > J_0, d_{J_i,i} < \varepsilon_1, \\ 0.5\left(\dfrac{1}{2^{\lfloor 0.5(J_{\max}+J_0) \rfloor}} + \dfrac{1}{2^{J_{\max}-1}}\right) & J_i = J_0, J_{i+1} > J_0, \\ \dfrac{1}{2^{J_{\max}-1}} & \text{else.} \end{cases} \quad (2.55)$$

For the WCU scheme, $h_i = 2^{-J_{\max}}$ is an appropriate value for all scalar conservation laws and Euler systems. However, for the AMWCU scheme, different strategies have to be applied as the above analysis. For the convenience of notation, TVBU denotes the reconstruction method for the WCU scheme, and TVBR and TVBC are applied to represent the reconstruction methods for the two types of problems with the AMWCU scheme.

2.5 Some remarks

The merits of the high-order upwind scheme, wavelet collocation methods, and adaptive algorithm are combined to design the present wavelet collocation upwind schemes. The "upwind" property stabilizes the solving by the implicit numerical viscosity and ensures the numerical solution converges to the physical solution correctly. Our numerical test shows that the orders of accuracy are consistent with the expected ones. The adaptive wavelet schemes with reconstruction process capture the shock waves by refining localized nodes. Therefore, the high-order adaptive wavelet schemes are very suitable for resolving the hyperbolic conservation laws whose solutions involve multiscale smooth structures and several discontinuities. Moreover, the adaptive wavelet collocation methods are more efficient compared with the Galerkin ones. The most attractive advantage is that the adaptive wavelet decomposition can approximate the function defined over an arbitrary geometry, which indicates that our schemes can be easily extended to complex domains in 2D or 3D problems.

3  Numerical experiments

In this section, several benchmark experiments are conducted by applying different schemes based on $N = 5$ and $N = 7$ wavelets. Two kinds of norms are defined to evaluate numerical errors as follows:

$$\|e\|_{l_\infty} = \max_k \left\{ |u_k^e - u_k^n| \right\}, \quad (3.1)$$

$$\|e\|_{l_2} = \left( \sum_k |u_k^e - u_k^n|^2 \Delta x \right)^{\frac{1}{2}}, \quad (3.2)$$

where $u_k^e$ is the exact solution, and $u_k^n$ is the numerical one.

## 3.1 Linear scalar equation in one dimension

Numerical examples are performed to the following 1D linear scalar equation:

$$u_t + au_x = 0, \quad x \in [-1,1]. \tag{3.3}$$

### 3.1.1 Accuracy verification

First, the order of accuracy of the proposed method in uniform node distribution is verified by refinement tests. A smooth initial condition is given by

$$u(x,0) = \sin(\pi x) \quad \text{periodic}. \tag{3.4}$$

The period boundary condition is treated by the extension method for periodic functions proposed in our previous work [71].

Numerical errors are computed at $t = 2$. Time steps that are small enough are selected so that the errors are dominated by the spatial discretization. The numerical errors and orders of accuracy for the two schemes are shown in Table 2. As mentioned in **Theorem 2.2**, the theoretical orders of accuracy of the schemes for $N = 5$ and $N = 7$ wavelets are fourth-order and sixth-order, respectively. The order of accuracy agrees with the expected one for the fourth-order scheme, and the numerical results show a higher order of accuracy than the theoretical one for the sixth-order scheme. The order of accuracy for the other methods, such as multimoment constrained FVM [12], WENO scheme [7], and spectral (finite) volume method [13], only provide the expected order of accuracy when the number of grid is large enough. Our schemes show a uniform convergence rate even for much coarser node distribution.

Table 2. Numerical errors and order of accuracy for 1D linear scalar equation

| Method | $N_1$ | $l_\infty$ error | $l_\infty$ order | $l_2$ error | $l_2$ order |
|---|---|---|---|---|---|
| $N = 5$ (4th order) | 16 | 2.99E−3 | — | 3.17E−3 | — |
| | 32 | 1.84E−4 | 4.02 | 1.90E−4 | 4.05 |
| | 64 | 1.15E−5 | 4.01 | 1.16E−5 | 4.03 |
| | 128 | 7.15E−7 | 4.00 | 7.21E−7 | 4.01 |
| | 256 | 4.47E−8 | 4.00 | 4.49E−8 | 4.01 |
| $N = 7$ (6th order) | 16 | 6.84E−5 | — | 7.05E−5 | — |
| | 32 | 1.02E−6 | 6.07 | 1.04E−6 | 6.08 |
| | 64 | 1.46E−8 | 6.12 | 1.48E−8 | 6.14 |
| | 128 | 1.71E−10 | 6.42 | 1.73E−10 | 6.42 |
| | 256 | 8.70E−13 | 7.61 | 8.75E−13 | 7.62 |

## 3.1.2 Advection of a square wave

To validate the ability of the schemes to capture a jump discontinuity, the experiment of advection of a square wave is carried out with the following initial condition:

$$u(x,0) = \begin{cases} 1 & -0.4 \leq x \leq 0.4, \\ 0 & \text{otherwise}, \end{cases} \quad \text{periodic}. \quad (3.5)$$

The numerical results at $t = 2$ and $t = 32$ are obtained by the WCU at the resolution level $J = 8$ without the limiting process in Figs. 5 and 6, respectively. Spurious oscillations pollute numerical solutions near the jump discontinuities for $N = 5$ and $N = 7$ schemes. The numerical results have a smaller oscillation amplitude and a thinner oscillatory region near the singularities for the higher-order wavelet upwind scheme. Fig. 6 also indicates that our schemes are stable for a long-term integration when solving linear scalar problems involving discontinuities without limiting.

The numerical results at $t = 2$ with the TVBU limiter at different uniform resolution levels are plotted in Fig. 7. The parameter M used in this test is listed in Table 3. Spurious oscillations are removed effectively, and the numerical solutions converge to the exact one gradually as $J$ increases. Thus, our WCU scheme with the TVBU limiter can capture discontinuities without visible spurious oscillations.

Table 3. Parameter M for different resolutions level $J_{max}$

| $J_{max}$ | 6 | 7 | 8 | 9 | 10 | 11 | 12 | 13 |
|---|---|---|---|---|---|---|---|---|
| M | 5 | 10 | 20 | 40 | 80 | 120 | 160 | 320 |

Then, the AMWCU scheme with the TVBR limiter is used to solve the problem for further validation. The numerical results are shown in Fig. 8. The solutions are non-oscillatory and approach the exact one with an increment in $J_{max}$. To evaluate the merits of the adaptive methods, the numerical results near the discontinuity of the classic fifth-order finite difference WENO scheme (WENO-5) proposed by Jiang and Shu [7] are compared with our adaptive algorithm, as illustrated in Fig. 9. Figs. 9(a) and (b) show that the proposed adaptive wavelet algorithm can capture steeper gradient than the WENO-5 scheme in much less nodes. The ratios of the node number of the AMWCU scheme to that of the WENO-5 scheme are about 10%, which shows considerable data compression. Sparse representation with higher order of accuracy is one of the main advantage of the proposed adaptive scheme, and the present adaptive scheme with the integral reconstruction can very efficiently capture discontinuities.

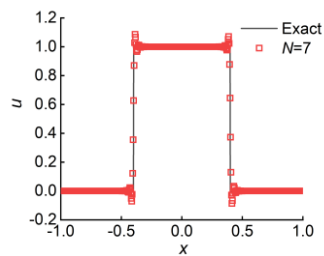
(a) $N=5$

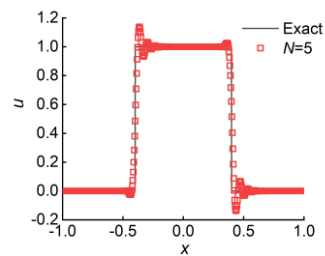
(b) $N=7$

Fig.5. Advection of a square wave by the WCU without limiting at $t = 2$s, $J = 8$

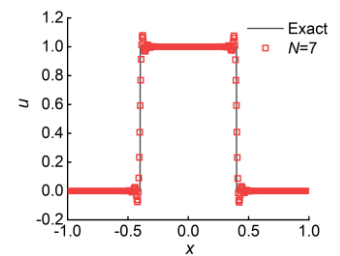
(a) $N=5$

(b) $N=7$

Fig.6. Advection of a square wave by the WCU without limiting at $t = 32$s, $J = 8$

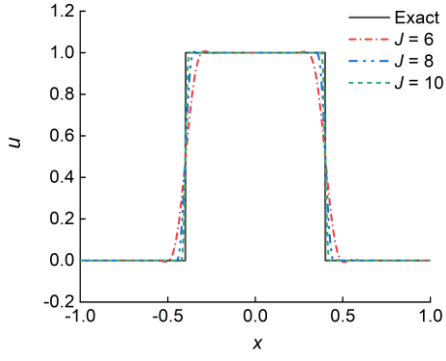
(a) N=5

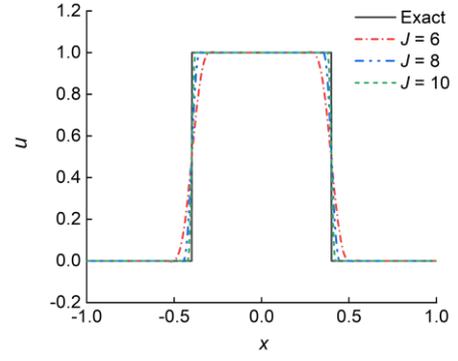
(b) N=7

Fig.7. Comparison of numerical solutions by the WCU with the TVBU at $t = 2$s

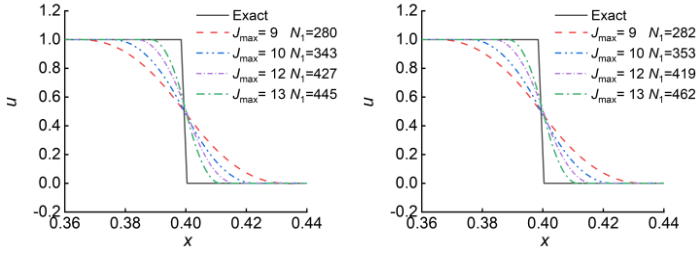

(a) N=5　　　　　　　　　(b) N=7

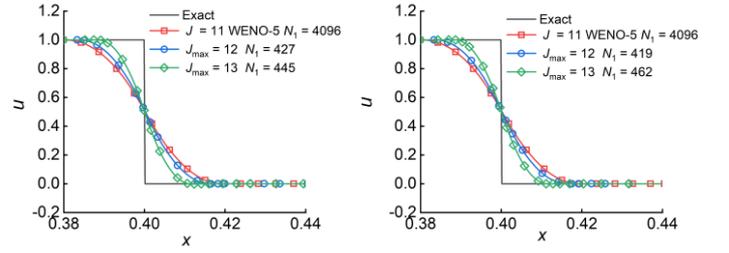

(a) N=5　　　　　　　　　(b) N=7

Fig.8. Comparison of numerical solutions by the AMWCU with the TVBR at $t = 2$s ($J_0 = 6$)

Fig.9. Comparison of the solutions between the AMWCU with the TVBR and WENO-5 ($J_0 = 6$)

3.1.3 Jiang and Shu's problem[7]

The example designed by Jiang and Shu [7] is used to investigate the capability of the wavelet schemes in recognizing smooth and discontinuous solutions. The initial condition contains a smooth but narrow combination of Gaussian, a square wave, a sharp triangle wave, and a half ellipse shown as

$$u(x,0) = \begin{cases} \frac{1}{6}\left(G(x,\beta,z-\delta)+G(x,\beta,z+\delta)+4G(x,\beta,z)\right) & -0.8 \leq x \leq -0.6, \\ 1 & -0.4 \leq x \leq -0.2, \\ 1-|10(x-0.1)| & 0 \leq x \leq 0.2, \\ \frac{1}{6}\left(F(x,\alpha,a-\delta)+F(x,\alpha,a+\delta)+4F(x,\alpha,a)\right) & 0.4 \leq x \leq 0.6, \\ 0 & \text{otherwise}, \end{cases} \quad (3.6)$$

$$G(x,\beta,z) = e^{-\beta(x-z)^2}, F(x,\alpha,a) = \sqrt{\max(1-\alpha^2(x-a)^2,0)}.$$

The constants are selected as $a = 0.5, z = -0.7, \delta = 0.005, \alpha = 10, \beta = \ln 2/(36\delta^2)$.

The numerical results at $t = 2$ and $t = 20$ are calculated by the WCU at the resolution level $J = 8$, as illustrated in Figs. 10 and 11, respectively. Fig. 10 shows that spurious oscillation appears in the local regions near the square wave and the half ellipse where the Lipschitz exponents are exactly zero. Moreover, the result of the $N = 7$ wavelet scheme exhibits better accuracy for smooth and discontinuous regions, especially for the jump discontinuity regions. Fig. 11 shows the $N = 7$ wavelet scheme produces a narrow oscillatory region and a lower overshoot in the vicinity of the square wave and a higher precision in smooth regions. The results of the $N = 5$ wavelet scheme reveals a slight distortion near the extreme point of the half ellipse.

Next, numerical experiments are carried out by employing the AMWCU scheme with the TVBR limiter. Fig. 12 shows the results of the two wavelet schemes at the basic resolution level $J_0 = 6$. The AMWCU scheme with the TVBR can efficiently capture the discontinuities and maintain high accuracy in smooth regions. A total of 6144 nodes are needed to improve the steepness of the square wave from $J = 10$ to $J = 12$ for the WCU scheme. With the AMWCU scheme, the goal can be realized by adding only 311 and 312 nodes for the $N = 5$ and $N = 7$ wavelet schemes, respectively, which shows the great ability of representing a function with different regularities.

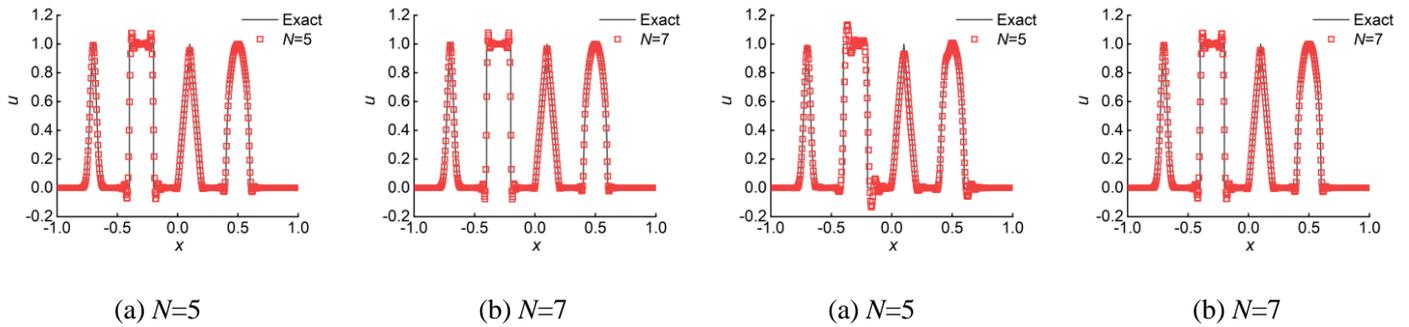

(a) $N=5$  (b) $N=7$  (a) $N=5$  (b) $N=7$

Fig.10. Numerical results by the WCU without limiting at $t = 2$, $J = 8$

Fig.11. Numerical results by the WCU without limiting at $t = 20$, $J = 8$

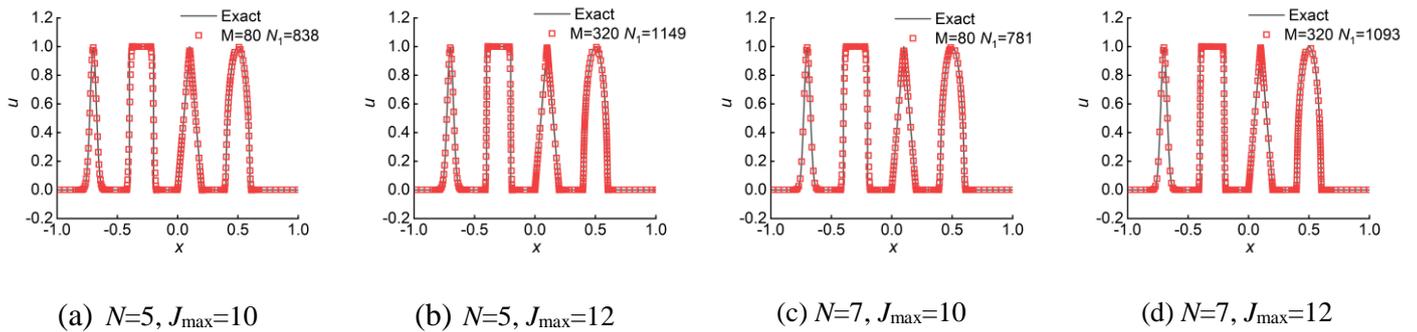

(a) $N=5$, $J_{max}=10$  (b) $N=5$, $J_{max}=12$  (c) $N=7$, $J_{max}=10$  (d) $N=7$, $J_{max}=12$

Fig.12. Comparison of numerical solutions of the AMWCU with the TVBR at $t = 2$ ($J_0 = 6$)

3.1.4 Advection of a mixing scale function

To show the ability of the present scheme in detecting different scale structures, the following initial condition with mixing low-frequency and high-frequency smooth sine functions and discontinuities are devised:

$$u(x,0) = \begin{cases} 0.5 & -0.8 \le x \le -0.6, \\ -2.5x - 0.5 & -0.6 \le x \le -0.4, \\ -2.5x & -0.4 \le x \le -0.2, \\ 0.5 & -0.2 \le x \le -0.1, \\ 0.5(1+\sin(40\pi x)) & -0.1 \le x \le 0.1, \\ 0.5 & 0.1 \le x \le 0.2, \\ 0.2\sin(5\pi(x-0.2)) & 0.2 \le x \le 0.6, \\ 0.5 & 0.6 \le x \le 0.8, \\ 0 & \text{otherwise} \end{cases} \quad (3.7)$$

The numerical results of the AMWCU are compared with those of the WENO-5 scheme, as illustrated in Fig. 13. Here the same display control parameters are set to improve the readability of the figure. The nodes of the AMWCU schemes are concentrated in discontinuous regions and high-frequency regions to depict different scales precisely, and the sixth-order scheme needs less nodes to obtain almost the same accuracy. However, for the WENO-5 scheme, a uniform node distribution with $N_1 = 2048$ is required to distinguish all the details of the solution. Fig. 14 shows the resolution level distribution of all nodes in the two schemes. The maximum resolution levels mark the jump discontinuity.

The solution of the test contains several different scales and discontinuities in a quite small domain. For such a case, the AMWCU schemes perform better in capturing all details with a highly efficient approach. Our AMWCU schemes possess an excellent ability to approximate complex solutions with high efficiency when addressing localized singular problems defined on large domains.

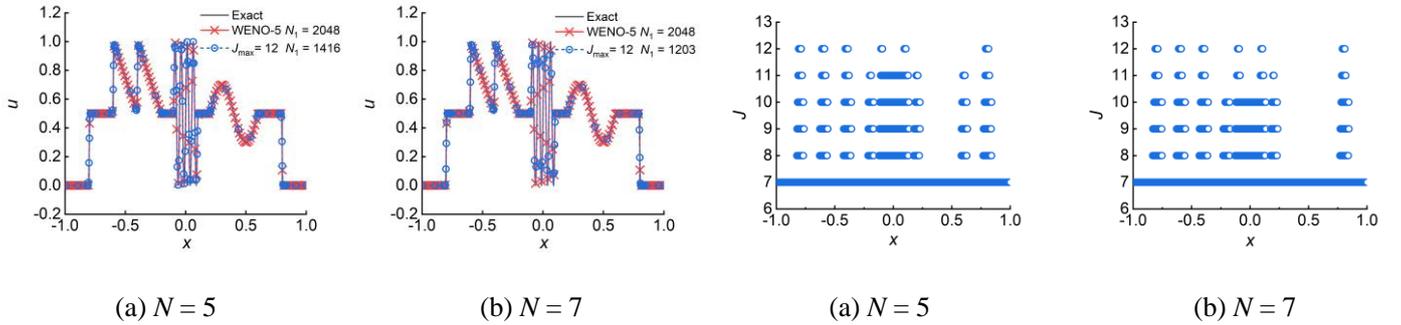

(a) $N = 5$    (b) $N = 7$    (a) $N = 5$    (b) $N = 7$

Fig.13. Comparison of numerical solutions by the AMWCU with the TVBR at $t = 2$ ($J_0=7$, M = 320)    Fig.14. Resolution levels by the AMWCU with the TVBR at $t = 2$ ($J_0 = 7$, M = 320)

3.2 Nonlinear inviscid Burger's equation in one dimension

The 1D nonlinear Burger's equation without viscosity is described as

$$\begin{aligned} u_t + \left(u^2/2\right)_x &= 0, & x \in [0,2], \\ u(x,0) &= 0.5 + \sin(\pi x), & \text{periodic.} \end{aligned} \quad (3.8)$$

This equation is a typical example that a shock wave develops as time increases with a sufficiently smooth initial condition owing to nonlinearity. First, the accuracy orders of the wavelet schemes are evaluated by node refinement tests. The exact solution is calculated by the semi-analytical method proposed by Harten et al. [4] that uses Newton–Raphson iterations to solve a characteristic relation. The exact solution at $t = 0.1$, which is still very smooth, is computed for comparison with the numerical ones. The errors and orders of accuracy in the $l_\infty$ and $l_2$ norms are shown in Table 4. The order of accuracy tends to the expected ones with an increase in node number.

Table 4. Numerical errors and orders of accuracy for 1D nonlinear inviscid Burger's equation

| Method | $N_1$ | $l_\infty$ error | $l_\infty$ order | $l_2$ error | $l_2$ order |
| --- | --- | --- | --- | --- | --- |
| $N = 5$ 4th order | 64 | 1.63E−5 | — | 8.06E−6 | — |
| | 128 | 1.24E−6 | 3.71 | 5.83E−7 | 3.79 |
| | 256 | 8.26E−8 | 3.91 | 3.84E−8 | 3.93 |
| | 512 | 5.29E−9 | 3.97 | 2.45E−9 | 3.97 |
| | 1024 | 3.34E−10 | 3.99 | 1.54E−10 | 3.99 |
| $N = 7$ 6th order | 64 | 2.74E−7 | — | 1.20E−7 | — |
| | 128 | 7.64E−9 | 5.17 | 3.07E−9 | 5.29 |
| | 256 | 1.46E−10 | 5.71 | 5.65E−11 | 5.77 |
| | 512 | 2.48E−12 | 5.87 | 9.54E−13 | 5.89 |
| | 1024 | 3.97E−14 | 5.97 | 2.22E−14 | 5.43 |

The numerical solutions are calculated up to $t = 1.5/\pi$ when a shock is developed. The results obtained by the WCU schemes with the TVBU at different resolution levels are plotted in Fig. 15. The solutions all converge to the exact ones without visible spurious oscillations as the resolution level increases. This finding demonstrates that our WCU schemes with the TVBU can work well in dealing with nonlinear hyperbolic problems and obtain a satisfactory physical solution.

To investigate influences of the TVBU limiters on the solution accuracy away from the discontinuity, the absolute errors over the interval of the two schemes at different resolution levels are displayed in Fig. 16. The expected order of accuracy for the $N = 5$ wavelet scheme with the TVBU limiter is retained away from the shock wave. Owing to the high accuracy of the $N = 7$ wavelet scheme, the errors in smooth regions sharply decline to machine error, inducing the desired accuracy order only at the low resolution levels.

Finally, a numerical experiment is conducted by using the AMWCU scheme with the TVBC. Fig. 17 illustrates the results at different maximum resolution levels. The numerical results of the adaptive scheme with the TVBC are slightly sharper than those of the WCU schemes with the basic resolution level, and the improvement of the accuracy seems restricted. This finding can be attributed to the reconstruction with large integral interval. Therefore, an optimized method to determine $h_i$ needs to be designed to remove the oscillations without smoothing the discontinuity too much.

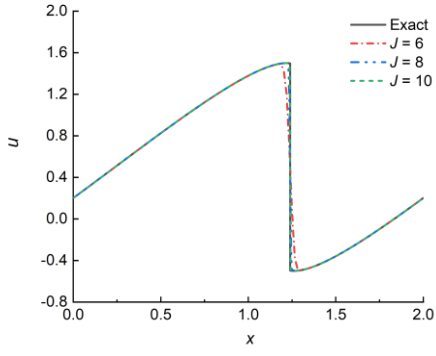
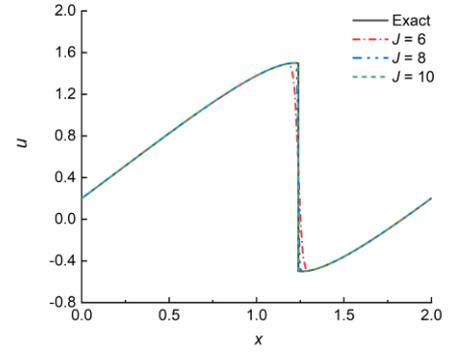

(a) $N = 5$　　　　　　　　　　　　　　　　　　(b) $N = 7$

Fig.15. Comparison of numerical solutions by the WCU with the TVBU at $t = 1.5/\pi$

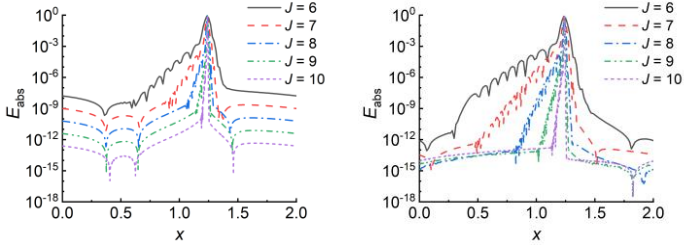
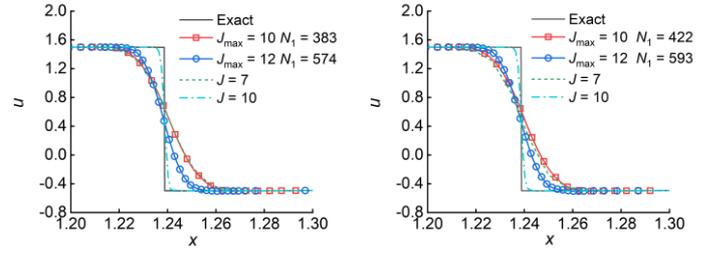

(a) $N = 5$　　　(b) $N = 7$　　　　　　(a) $N = 5$　　　(b) $N = 7$

Fig.16. Comparison of numerical errors by the WCU with the TVBU at $t = 1.5/\pi$　　　Fig.17. Comparison of numerical solutions by the AMWCU with the TVBC at $t = 1.5/\pi$ ($J_0 = 7$)

### 3.3 Euler system in one dimension

The 1D Euler equations of gas dynamics for a polytropic gas in the conservative form are written as

$$\begin{aligned}&\mathbf{U}_t + \mathbf{f}(\mathbf{U})_x = 0, \\ &\mathbf{U} = (\rho, \rho u, E)^T, \\ &\mathbf{f}(\mathbf{U}) = (\rho u, \rho u^2 + p, u(E + p))^T,\end{aligned} \qquad (3.9)$$

where $\rho$, $u$, $p$, and $E$ are the density, velocity, pressure, and total energy, respectively. For the perfect gas, the state equation is

$$E = \frac{p}{\gamma - 1} + \frac{1}{2}\rho u^2, \qquad (3.10)$$

where $\gamma = 1.4$.

#### 3.3.1 Advection of density perturbation

To check the order of accuracy of the proposed schemes for the nonlinear systems, the propagation of density perturbation test devised by Qiu and Shu is computed [72]. The initial condition is given as $\rho(x, 0) = 1 + 0.2\sin(\pi x)$, $u(x, 0) = 1$, $p(x, 0) = 1$ over $[0, 2]$ with periodic boundary condition.

The resolution level is refined from 3 to 7. The numerical errors at $t = 2$ evaluated by the density are shown in Table 5. The convergence rates are all higher than the theoretical ones. For the $N = 5$ wavelet scheme, the convergence rate decreases to the expected one as the resolution increases. However, a super convergence is revealed for the $N = 7$ wavelet scheme with an increment in the resolution level.

Table 5. Numerical errors and orders of accuracy for 1D Euler system

| Method | $N_1$ | $L_\infty$ error | $L_\infty$ order | $L_2$ error | $L_2$ order |
|---|---|---|---|---|---|
| $N=5$ (4th order) | 16 | 8.86E−4 | — | 9.19E−4 | — |
| | 32 | 4.28E−5 | 4.37 | 4.37E−5 | 4.39 |
| | 64 | 2.39E−6 | 4.16 | 2.43E−6 | 4.17 |
| | 128 | 1.45E−7 | 4.05 | 1.46E−7 | 4.06 |
| | 256 | 8.97E−9 | 4.01 | 9.00E−9 | 4.02 |
| $N=7$ (6th order) | 16 | 2.49E−5 | — | 2.52E−5 | — |
| | 32 | 2.63E−7 | 6.56 | 2.67E−7 | 6.56 |
| | 64 | 3.21E−9 | 6.36 | 3.25E−9 | 6.36 |
| | 128 | 3.59E−11 | 6.47 | 3.61E−11 | 6.49 |
| | 256 | 3.69E−13 | 6.61 | 2.99E−13 | 6.91 |

3.3.2 Sod's problem

This experiment is a well-known standard test as a Riemann problem proposed by Sod [73] with the following initial data:

$$(\rho_0, u_0, p_0) = \begin{cases} (1,\ 0,\ 1) & 0 \leq x \leq 0.5 \\ (0.125,\ 0,\ 0.1) & \text{otherwise} \end{cases} \quad (0 \leq x \leq 1). \tag{3.11}$$

First, the numerical results at $t = 0.2$ are calculated by applying the WCU scheme without limiting, as shown in Fig. 18. The wavelet scheme without a reconstruction process can locate the shock discontinuity and contact discontinuity positions accurately, and only local spurious oscillations arise. Then, the TVBU limiter is added to remove the numerical oscillations. Fig. 19 displays the density obtained by the WCU schemes with the TVBU. Spurious oscillations near the discontinuities are eliminated successfully.

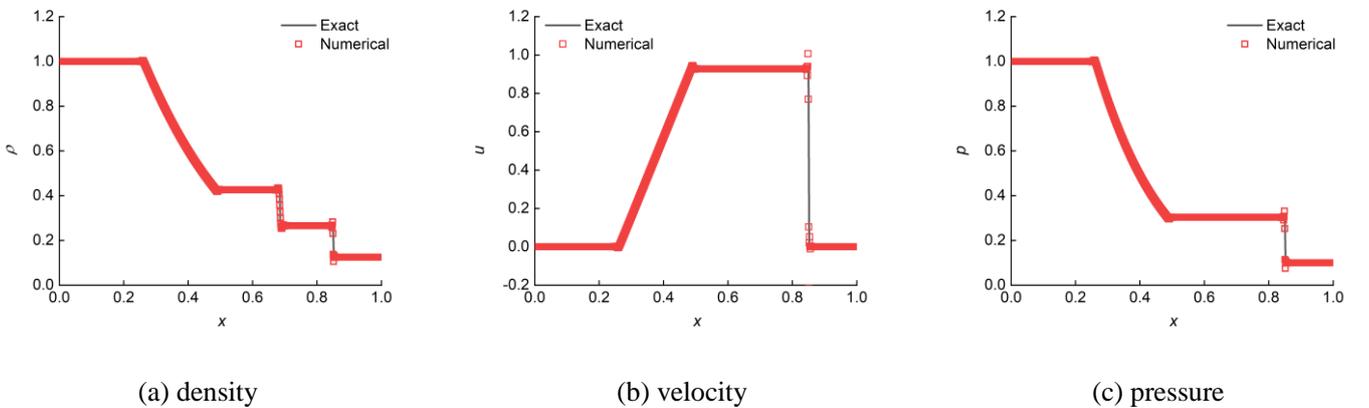

(a) density      (b) velocity      (c) pressure

Fig.18. Numerical solutions by the WCU without limiting at $t = 0.2$ ($N = 5$, $J = 10$)

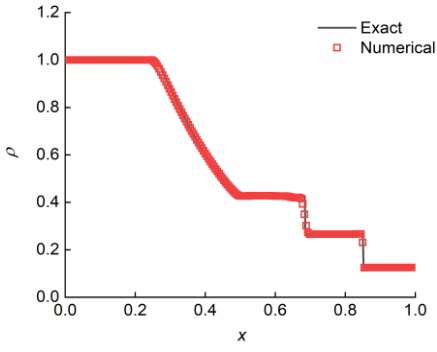
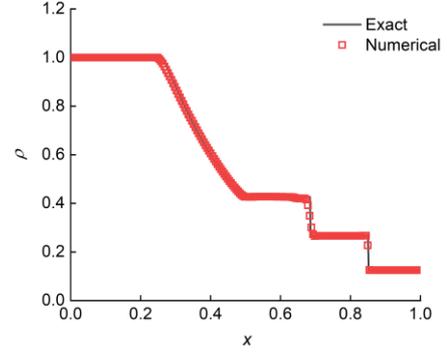

(a) $N = 5$  (b) $N = 7$

Fig.19. Comparison of density by the WCU with the TVBU at $t = 0.2$ ($J = 10$, M=40)

The performance of the AMWCU with the TVBR on such a Riemann problem is further explored. Figs. 20 and 21 compare the density obtained with schemes of different orders and maximum resolution levels. The solutions with $J_{max} = 12$ behave much steeper in the vicinity of the jump discontinuities, and the higher-order wavelet scheme manifests a better ability to suppress the Gibbs phenomenon at the lower maximum resolution level.

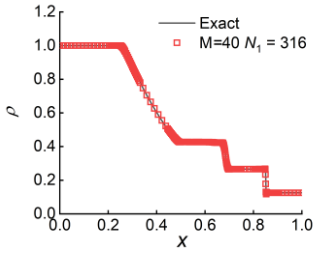
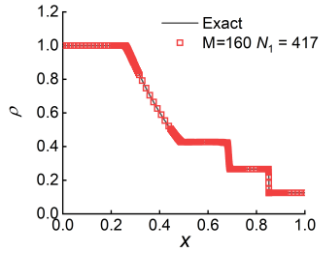
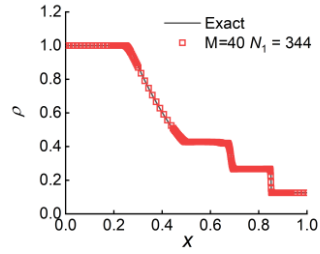
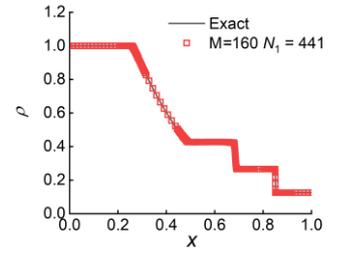

(a) $J_{max}=10$    (a) $J_{max}=12$    (a) $J_{max}=10$    (b) $J_{max}=12$

Fig.20. Numerical solutions by the AMWCU With the TVBR at $t = 0.2$ ($N = 5$, $J_0 = 6$)

Fig.21. Numerical solutions by the AMWCU With the TVBR at $t = 0.2$ ($N = 7$, $J_0 = 6$)

3.3.3 Lax's problem

Another famous Riemann problem is devised by Lax [74] with the following initial condition:

$$(\rho_0, u_0, p_0) = \begin{cases} (0.445,\ 0.698,\ 3.528) & 0 \leq x \leq 0.5 \\ (0.5,\ 0,\ 0.571) & \text{otherwise} \end{cases} \quad (0 \leq x \leq 1). \qquad (3.12)$$

The numerical results at $t = 0.13$ are computed based on the WCU scheme without the reconstruction. The distributions of three variables are illustrated in Fig. 22. Numerical oscillations are confined to the local region of the contact and shock discontinuities, and the density reveals the severest oscillatory characteristic. The original WCU scheme can stably solve Lax's problem and accurately recognize all types of discontinuities.

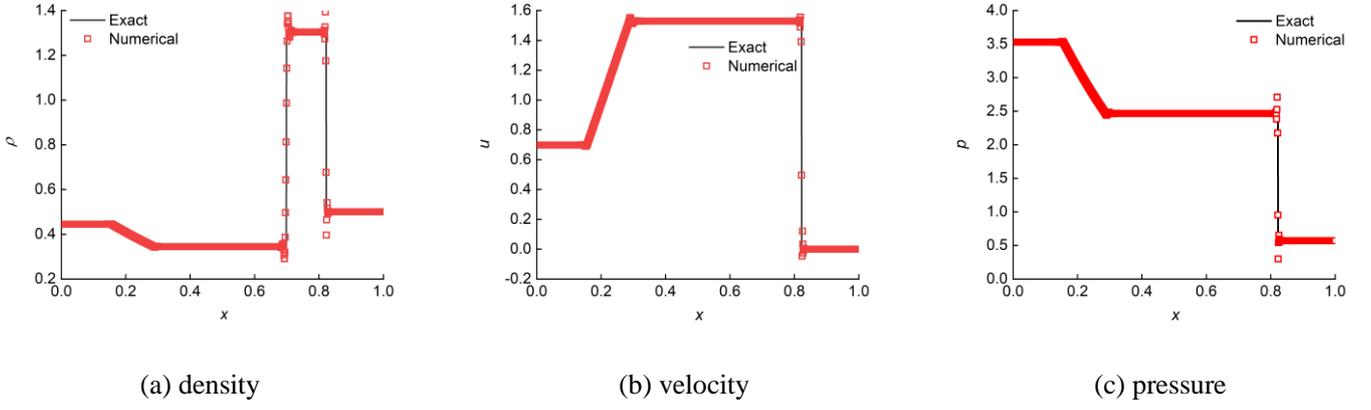

(a) density  (b) velocity  (c) pressure

Fig.22. Numerical solutions by the WCU without limiting at $t = 0.13$ ($N = 5$, $J = 10$)

Next, experiments are conducted by employing the AMWCU with the TVBR. Figs. 23 and 24 display the numerical results of density with different wavelet schemes and maximum resolution levels. Jiang and Shu [7] pointed out that this test is so tough for high-order non-characteristic-based schemes that oscillations could easily appear. In Figs. 23 and 24, the numerical solutions of the $N = 7$ wavelet scheme are ENO, which shows the extraordinary capacity of the AMWCU scheme in solving the hyperbolic conservation laws with discontinuous solutions. For the $N = 5$ scheme, only small acceptable oscillations emerge in the very thin region of the shock wave front. Thus, the high-order scheme with the TVBR performs well in removing the oscillations near the shock discontinuities. Comparing the rarefaction wave of the solutions in different order schemes reveals that the higher-order scheme can accurately approximate the rarefaction wave only at the basic resolution level. Hence, it requires less nodes but shows better performance in solving the current problem.

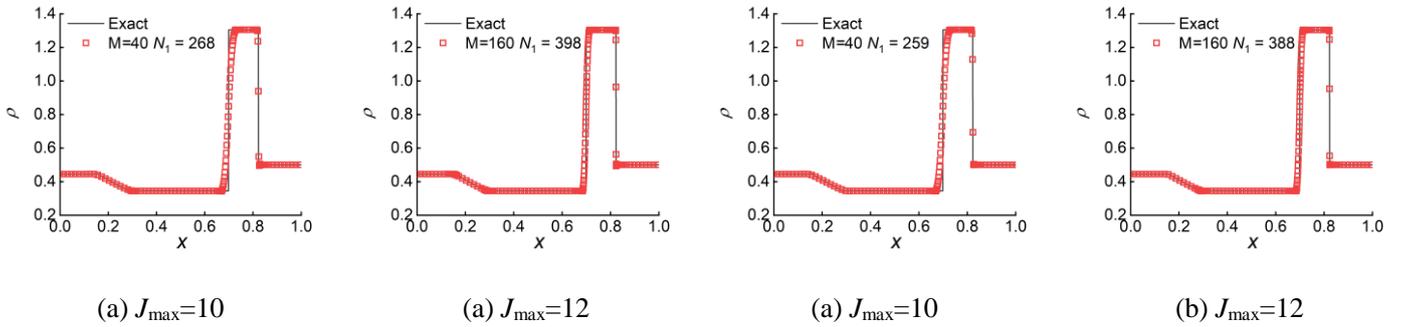

(a) $J_{max}=10$  (a) $J_{max}=12$  (a) $J_{max}=10$  (b) $J_{max}=12$

Fig.23. Numerical solutions by the AMWCU With the TVBR at $t = 0.13$ ($N = 5$, $J_0 = 7$)

Fig.24. Numerical solutions by the AMWCU With the TVBR at $t = 0.13$ ($N = 7$, $J_0 = 7$)

3.3.4 Shock–turbulence interaction

High-order schemes demonstrate their merits when the solution contains shocks and complex smooth region structures [20]. Shock interaction with entropy waves is such a typical example proposed by Shu and Osher [75]. The problem of a moving Mach = 3 shock interacting with sine waves in density is defined by the following initial data:

$$(\rho_0, u_0, p_0) = \begin{cases} (3.857148,\ 2.629369,\ 10.333333) & -5 \leq x < -4 \\ (1+0.2\sin 5x,\ 0,\ 1) & \text{otherwise} \end{cases} \quad (-5 \leq x \leq 5) \quad (3.13)$$

In the present case, the reference "exact" solution, which is a converged solution computed by the WENO-5 scheme with 2560 grid points, is obtained. The calculated density $\rho$ based on the WCU scheme without limiting against the reference solution is plotted at $t = 1.8$ in Fig. 25. The wavelet schemes can capture shocks and complex smooth structures at different scales. The problem is solved applying the AMWCU schemes with the TVBR. Comparison of the numerical results and the reference solution is illustrated in Figs. 26 and 27. Nodes at high resolution level gather in the vicinity of discontinuities and in high frequency regions, and the AMWCU schemes with the reconstruction can remove numerical oscillations effectively. Thus, the adaptive wavelet schemes perform quite well in distinguishing complex structures with shocks.

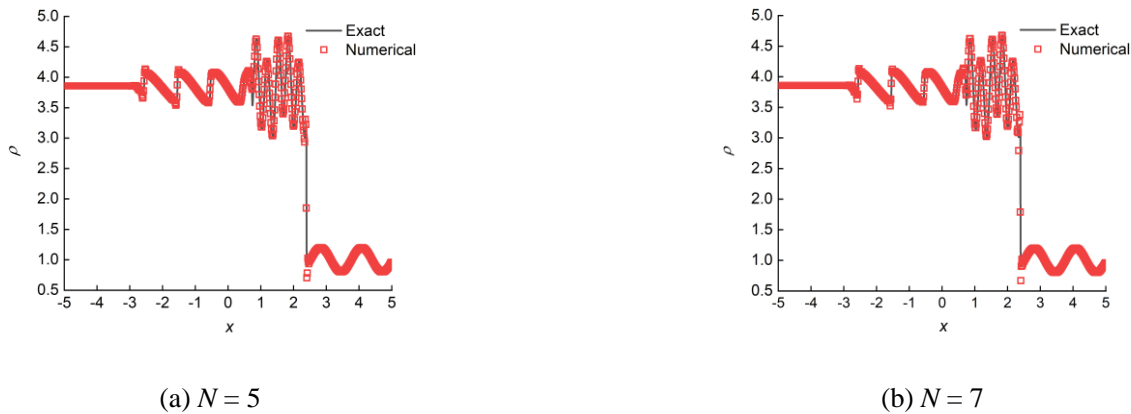

(a) $N = 5$  (b) $N = 7$

Fig.25. Numerical solutions by the WCU without limiting at $t = 1.8$ ($J = 6$)

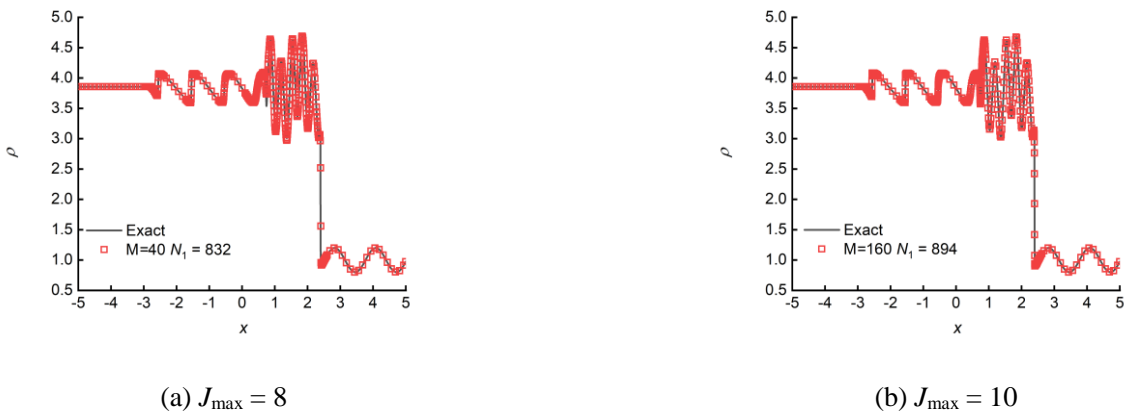

(a) $J_{max} = 8$  (b) $J_{max} = 10$

Fig.26. Numerical solutions by the AMWCU With the TVBR at $t = 1.8$ ($N = 5$, $J_0 = 4$)

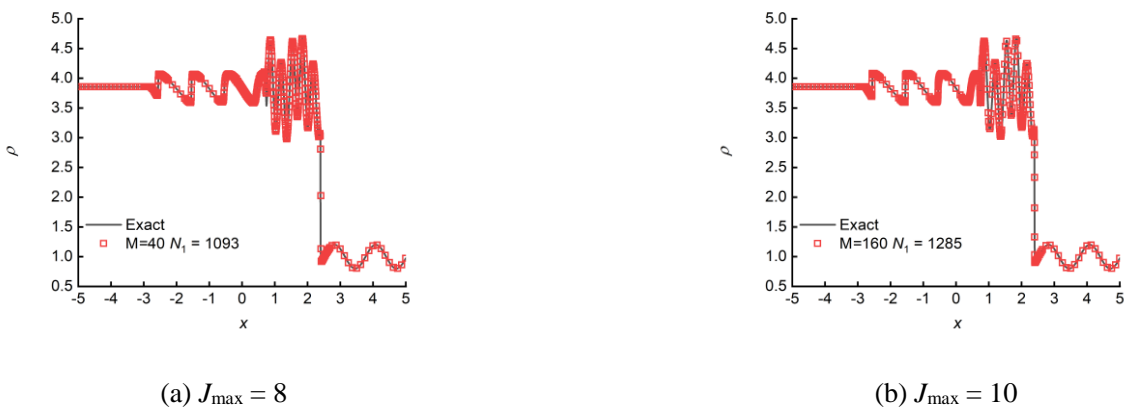

(a) $J_{max} = 8$  (b) $J_{max} = 10$

Fig.27. Numerical solutions by the AMWCU With the TVBR at $t = 1.8$ ($N = 7$, $J_0 = 4$)

### 3.3.5 Two interacting blast waves

Finally, the blast wave problem [76] involving complex interactions of strong shock waves, rarefaction waves, and contact discontinuities with one another is considered. The initial condition of the problem with substantial large pressure difference is described as follows:

$$(\rho_0, u_0, p_0) = \begin{cases} (1, \ 0, \ 1000) & 0 \leq x < 0.1 \\ (1, \ 0, \ 0.01) & 0.1 \leq x < 0.9 \\ (1, \ 0, \ 100) & \text{otherwise} \end{cases} \quad (0 \leq x \leq 1). \tag{3.14}$$

A reflective boundary condition is imposed at both ends, which can be seen in [4, 76].

This test is quite challenging because the complicated interactions of all types of waves are very strong and instantaneous. Here the WENO-5 scheme is still employed with sufficiently fine nodes to compute the reference "exact" solution. To validate the ability of the wavelet schemes in solving the present problem, the numerical results are computed up to $t = 0.038$ by applying the WCU schemes with the TVBR. The numerical results of density compared with the reference solution are shown in Figs. 28 and 29. The higher-order scheme behaves more accurately in tracing the peak and the trough. Moreover, the numerical solutions of the wavelet schemes converge to the reference "exact" one when increasing the resolution level.

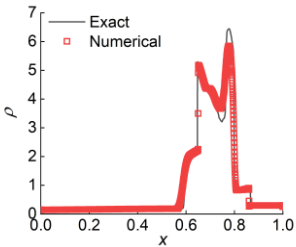 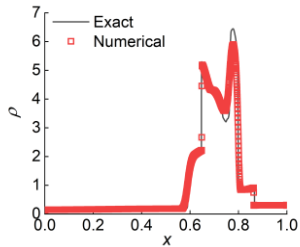 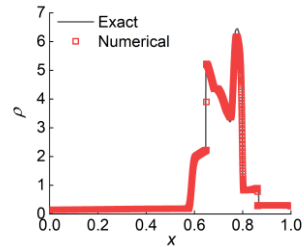 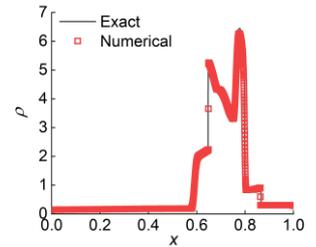

(a) $N = 5$    (b) $N = 7$    (a) $N = 5$    (b) $N = 7$

Fig.28. Numerical solutions by the WCU with the TVBU at $t = 0.038$ ($J = 12$, M = 3200)    Fig.29. Numerical solutions by the WCU with the TVBU at $t = 0.038$ ($J = 14$, M = 8000)

Next, the numerical results obtained by the adaptive wavelet schemes are presented. Figs. 30 and 31 illustrate the computed density obtained by the AMWCU schemes with different orders against the "exact" solution. The results of the AMWCU coincide with the solutions calculated by the WCU with the resolution level $J_{\max}$. The numerical results suggest that the adaptive wavelet schemes can be employed to resolve complex wave interaction problems efficiently.

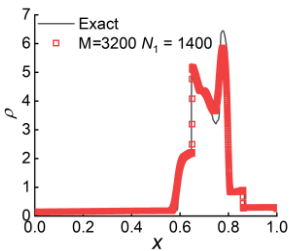 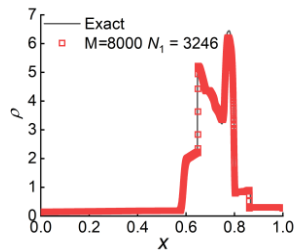 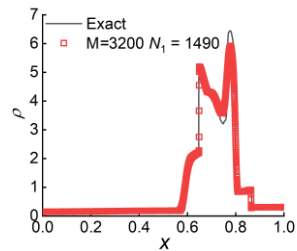 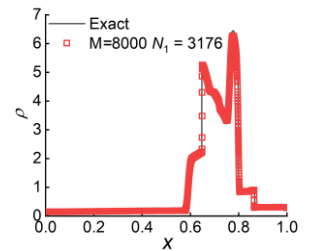

(a) $J_{\max} = 12$    (b) $J_{\max} = 14$    (a) $J_{\max} = 12$    (b) $J_{\max} = 14$

Fig.30. Numerical solutions by the AMWCU With the TVBR at $t = 0.038$ ($N = 5$, $J_0 = 8$)    Fig.31. Numerical solutions by the AMWCU With the TVBR at $t = 0.038$ ($N = 7$, $J_0 = 8$)

# 4  Conclusion

High-order adaptive multiresolution wavelet collocation upwind schemes for hyperbolic conservation laws based on asymmetrical interpolation wavelets are proposed. The upwind property can be achieved by using a couple of asymmetrical scaling functions. Thus, the flux splitting technique can be easily embedded into the wavelet schemes for 1D Euler system. Nonlinear terms in hyperbolic conservation laws can be addressed conveniently with high efficiency due to the interpolation property of the scaling functions.

Many numerical experiments in one dimension are performed to validate advantages of the present schemes. The orders of accuracy for the linear and nonlinear scalar equations and the Euler system agree with the expected orders. The numerical solutions for the shock turbulence interaction problem indicate that the proposed adaptive wavelet collocation upwind schemes with reconstruction process can stably resolve hyperbolic conservation laws and efficiently obtain high-order accuracy in smooth regions and steep discontinuity transitions without visible oscillations. Furthermore, our schemes can solve complex wave interaction problems, and the basic idea presented in this paper can be naturally extended to 2D or 3D problems.

## Acknowledgement

This research was supported by grants from the National Natural Science Foundation of China (11925204).